\documentclass[12pt]{article}
\usepackage{graphicx,psfrag,amsmath,amssymb}
\newcommand{\qed}{\nobreak \ifvmode \relax \else
\ifdim\lastskip<1.5em \hskip-\lastskip
\hskip1.5em plus0em minus0.5em \fi \nobreak
\vrule height0.75em width0.5em depth0.25em\fi}
\begin{document}

\title{Ramsey Functions for Generalized Progressions}

\author{Mano Vikash Janardhanan \& Sujith Vijay \\ 
IISER Thiruvananthapuram \\
\tt{manovikash@iisertvm.ac.in, sujith@iisertvm.ac.in}}

\date{}

\maketitle

\vskip 20pt

\centerline{\bf Abstract}

\vskip 5pt

\noindent

Given positive integers $n$ and $k$, a $k$-term semi-progression of scope 
$m$ is a sequence $(x_1,x_2,...,x_k)$ such that $x_{j+1} - x_j \in 
\{d,2d,\ldots,md\}, 1 \le j \le k-1$, for some positive integer $d$. Thus 
an arithmetic progression is a semi-progression of scope $1$. Let $S_m(k)$ 
denote the least integer for which every coloring of $\{1,2,...,S_m(k)\}$ 
yields a monochromatic $k$-term semi-progression of scope $m$. We obtain 
an exponential lower bound on $S_m(k)$ for all $m=O(1)$. Our approach also 
yields a marginal improvement on the best known lower bound for the 
analogous Ramsey function for quasi-progressions, which are sequences 
whose successive differences lie in a small interval.

\thispagestyle{empty}
\baselineskip=15pt
\vskip 30pt

\section*{\normalsize 1. Introduction}

In 1927, B.L. van der Waerden \cite{Wae27} proved that given positive 
integers $r$ and $k$, there exists an integer $W(r,k)$ such that any 
$r$-coloring of $\{1,2,\ldots,W(r,k)\}$ yields a monochromatic $k$-term 
arithmetic progression. Even after nine decades, the gap between the lower 
and upper bounds is enormous, with the best known lower bound of the order 
of $r^k$, whereas the best known upper bound is a five-times iterated 
tower of exponents (see \cite{Gow01}). Analogues of the Van der Waerden 
threshold $W(r,k)$ have been studied for many variants of arithmetic 
progressions, including semi-progressions and quasi-progressions (see 
\cite{Lan04}). \\

Given positive integers $m$ and $k$, a $k$-term semi-progression of scope 
$m$ is a sequence $(x_1,x_2,\ldots,x_k)$ such that for some positive 
integer $d$, $x_{j+1} - x_j \in \{d,2d,\ldots,md\}$. The integer $d$ is 
called the {\em {low-difference}} of the semi-progression. We define 
$S_m(k)$ as the least integer for which any $2$-coloring of 
$\{1,2,\ldots,S_m(k)\}$ yields a monochromatic $k$-term semi-progression 
of scope $m$. Note that $S_m(k) \le W(k)$ with equality if $m=1$.

\section*{\normalsize 2. An Exponential Lower Bound for $S_m(k)$}

Landman $\cite{Lan98}$ showed that $S_m(k) \ge (2k^2/m)(1+o(1))$. We 
improve this to an exponential lower bound for all $m=O(1)$. \\

{\bf {Theorem}} $S_m(k) > \alpha^k$ where $\alpha = \alpha(m) = 
\sqrt{2^m/(2^m-1)}$ \\

{\bf {Proof}} Let $f(N,k,m)$ denote the number of $2$-colorings of $[1,N]$ 
with a monochromatic $k$-term semiprogression of scope $m$. (In the 
remainder of the proof, we only consider $k$-term semi-progressions of 
scope $m$.) Note that $S_m(k)$ is the least integer $N$ such that 
$f(N,k,m)=2^N$. We derive an upper bound on $f(N,k,m)$ as follows. \\

Given a semi-progression $P=\{a_1,a_2,\ldots,a_k\}$ of low-difference $d$, 
we define the conjugate vector of $P$ as $(u_1,u_2,\ldots,u_{k-1})$ where 
$u_i=(a_{i+1}-a_i-d)/d$. Likewise, the frequency vector of $P$ is defined 
as $(v_0,v_1,\ldots,v_{m-1})$ where $v_j$ is the number of times $j$ 
occurs in the conjugate vector of $P$. Finally, the weight of $P$, denoted 
$w(P)$ is defined as $u_1+u_2+ \ldots+u_{k-1}$. \\

Given a coloring $\chi$, we define the {\em {$(a,d)$-primary 
semi-progression}} of $\chi$ as the semi-progression $P$ whose conjugate 
vector is lexicographically least among the conjugate vectors of all 
semi-progressions (with first term $a$ and low-difference $d$) that are 
monochromatic under $\chi$. Let $P=\{a_1,a_2,\ldots,a_k\}$ be a 
semi-progression with first term $a_1=a$ and low-difference $d$. We will 
give an upper bound for the number of colorings $\chi$ such that $P$ is 
the {\em {$(a,d)$-primary semi-progression}} of $\chi$. \\

Since $P$ is monochromatic, all elements of $P$ have the same color under 
$\chi$. Furthermore, if $(u_1,u_2,\ldots,u_{k-1})$ is the conjugate vector 
of $P$, it follows from the fact that $P$ is the $(a,d)$-primary 
semi-progression of $\chi$ that $w(P)$ elements in the arithmetic 
progression $\{a,a+d,\ldots,a+m(k-1)d\}$ must be of the color different 
from the color of the elements of $P$. For example, let $a=17,d=5,m=3,k=6$ 
and $P=\{17, 32, 42, 47, 62, 72\}$ with conjugate vector $(2,1,0,2,1)$. If 
the two colors are red and blue, and the elements of $P$ are all red, then 
$22,27,37,52,57$ and $67$ must all be blue. Indeed, if $57$ is red, then 
the semi-progression $P'=\{17,32,42,47,57,62\}$ would have a 
lexicographically lower conjugate vector $(2,1,0,1,0)$. Thus there are at 
most $2^{N-11}$ colorings of $[1,N]$ whose $(a,d)$-primary 
semi-progression is $P$. \\

Let ${\bf{\rho}}=(1,1,\ldots,1) \in {\mathbb Z}^m$ and 
${\bf{\mu}}=(0,1,\ldots,m-1)$. Clearly, $w(P)= \sum_{j=0}^{m-1} j v_j = 
<{\bf{\mu}},{\bf{v}}>$ where ${\bf{v}}$ is the frequency vector of $P$. 
Note that there are at most $N^2/(k-1)$ choices for the pair $(a,d)$. 
We say that two progressions $P_1$ and $P_2$ with the same $a$ and $d$ are 
equivalent if they have the same frequency vector.
Note that for any $a$ and $d$, there are at most $$M(P) = 
\frac{(v_0 + v_1 + \cdots + v_m)!}{v_0 ! v_1 ! \cdots v_{m-1} !}$$ 
semi-progressions with the same frequency vector 
$(v_0,v_1,\ldots,v_{m-1})$ as $P$. Adding over all the equivalence classes 
of semi-progressions, we obtain 
$$f(N,k,m) \le \frac{N^2 2^{N-k+1}}{k-1} 
\sum_{w(P)=0}^{(m-1)(k-1)} M(P) 2^{-w(P)}$$ It follows from the 
multinomial theorem that 
$$f(N,k,m) \le \frac{N^2 2^N}{k-1} {\left( \frac{1}{2} + \frac{1}{2^2} + 
\cdots + \frac{1}{2^m} \right) }^k$$ Thus $f(N,k,m) < 2^N$ for $N = 
\alpha_m^k$ where $\alpha_m = \sqrt{2^m/(2^m-1)}$. This completes the 
proof. \qed

\section*{\normalsize 3. Exponential Lower Bounds for $Q_n(r,k)$}

We now apply the same technique to quasi-progressions. A $k$-term 
quasi-progression of low difference $d$ and diameter $n$ is a sequence 
$(a_1,a_2,\ldots,a_k)$ such that $d \le a_{j+1} - a_j \le d+n, 1 \le j \le 
k-1$. Let $Q_n(r,k)$ denote the least positive integer such that any 
$r$-coloring of $\{1,2,\ldots,Q_n(r,k)\}$ yields a monochromatic $k$-term 
quasi-progression of diameter $n$. It is known (see \cite{Vij10}) that 
$Q_1(2,k) > {\beta}^k$ where $\beta=1.08226...$ is the smallest positive 
real root of the equation 
$$y^{24}+8y^{20}-112y^{16}-128y^{12}+1792y^8+1024y^4-4096=0$$ and that 
$Q_n(k) = O(k^2)$ for $n > k/2$ (see \cite{Job11}). We apply the 
techniques of the previous section to obtain lower bounds on $Q_n(r,k)$. 
Let $g(r,N,k,n)$ denote the number of $r$-colorings of $[1,N]$ with a 
monochromatic $k$-term semiprogression of diameter $n$. Note that 
$Q_n(r,k)$ is the least positive integer $N$ such that $g(r,N,k,n)=2^N$. 
We first discuss the simplest non-trivial case, namely $r=2$ and $n=1$. \\

We define the conjugate vector of a quasiprogression 
$Q=\{a_1,a_2,\ldots,a_k\}$ of low-difference $d$ as 
$(u_1,u_2,\ldots,u_{k-1})$ where $u_i=a_{i+1}-a_i-d$. Given a coloring 
$\chi$, we define the {\em {$(a,d)$-primary quasi-progression}} of $\chi$ 
as the quasi-progression $Q$ whose conjugate vector is lexicographically 
least among the conjugate vectors of all quasi-progressions (with first 
term $a$ and low-difference $d$) that are monochromatic under $\chi$. Let 
$Q=\{a_1,a_2,\ldots,a_k\}$ be a quasi-progression with first term $a_1=a$ 
and low-difference $d$. We give an upper bound for the number of colorings 
$\chi$ such that $Q$ is the {\em {$(a,d)$-primary quasi-progression}} of 
$\chi$. \\

Since $Q$ is monochromatic, all elements of $Q$ have the same color under 
$\chi$, say red. Let $(u_1,u_2,\ldots,u_{k-1})$ be the conjugate vector of 
$Q$. Observe that if $u_j=1$ and $u_{j+1}=0$ for some $j$, so that $a_j, 
a_j+d+1$ and $a_j+2d+1$ are elements of $Q$, and therefore red, it follows 
that the color of $a_j+d$ is different from red (say blue), as $(P \; \cup 
\; \{a_j+d\}) \; \setminus \; \{a_j+d+1\}$ has a lexicographically lower 
conjugate vector. We define the weight of $Q$, denoted $w(Q)$, as the 
sum of the last element of the conjugate vector of $Q$, and the number of 
occurrences of the string ``10" in the conjugate vector of $Q$. 
Note that in view of the above observation, the color of $w(Q)$ integers 
in the set $\{a, a+d, a+d+1, \ldots, a+(k-1)d, \ldots, a+(k-1)(d+1)\}$ can 
be inferred to be blue. \\

We now derive an upper bound on $g(2,N,k,1)$. There are $N^2/(k-1)$ 
choices for the pair $(a,d)$. Of the $2^{k-1}$ possible conjugate vectors 
for a quasi-progression with first term $a$ and common difference $d$, let 
$w_{\ell}$ be the number of conjugate vectors of weight $\ell$. Let 
$$S_{t}=\sum_{\ell=0}^{\lceil t/2 \rceil} w_{\ell} 2^{- \ell}$$ denote 
the weighted sum of all such vectors of length $t$. Clearly, $S_t=S_{t,0} 
+ S_{t,1}$ where $S_{t,0}$ and $S_{t,1}$ denote the weighted sum of 
conjugate vectors that begin with $0$ and $1$ respectively, with 
$S_{1,0}=1$ and $S_{1,1}=1/2$. It is easy to see that $A {[S_{t-1,0} \; 
S_{t-1,1}]}^T = {[S_{t,0} \; S_{t,1}]}^T$ where $$A=\left[ 
\begin{array}{cc} 1 & 1 \\ 1/2 & 1 \end{array} \right]$$ Since 
$\lambda_{max}(A)=1+\frac{1}{\sqrt{2}}$, we get $$g(2,N,k,1) < \frac{N^2 
2^{N-k+1} \left[ \left(1+\frac{1}{\sqrt{2}} \right)^k + 
\left(1-\frac{1}{\sqrt{2}} 
\right)^k \right] }{2(k-1)}$$ Thus $g(2,N,k,1) < 
2^N$ for $N = \beta_{2,1}^k$ where $\beta_{2,1} = 1.08239...\,$ is the 
smallest positive real root of the equation $y^4-8y^2+8=0$. It follows 
that $Q_1(2,k) > \beta^k_{2,1}$ yielding a marginal improvement over the 
lower bound in \cite{Vij10}. \\

In general, since there are $r^N$ $r$-colorings of $[1,N]$ and at most 
$N^2(n+1)^{k-1}$ $k$-term quasi-progressions of diameter $n$, a lower 
bound of the form $Q_n(r,k) \ge {(\sqrt{r/(n+1)})}^k$ follows immediately 
from the linearity of expectation. However, this bound is only useful when 
$n \le r-2$. Generalising the approach outlined earlier, we represent the 
conjugate vector of $Q$ as an $r$-ary string, and define the weight $w(Q)$ 
as the sum of the last element of the conjugate vector of $Q$, and the 
number of occurrences of strings of length two of the form ``{\em{xy}}", 
counted with multiplicity $m(x,y)=\min(x,n-y)$. (Note that $m(x,y)$ 
denotes the number of conjugate vectors that are lexicographically lower 
than the given vector and correspond to quasi-progressions that differ 
from $Q$ in exactly one element.) As before, let $S_{t,j}$ 
denote the weighted sum of of conjugate vectors of length $t$ beginning 
with $j, 0 \le j \le n$, with $S_{1,j}=\alpha^j$ for 
all $j$ where $\alpha=1-\frac{1}{r}$. Then $A [S_{t,0} \; \cdots \; S_{t,n}]^T 
= [S_{t+1,0} \; \cdots \; S_{t+1,n}]^T$ where $$A_{r,n}=\left[ 
\begin{array}{ccccc} 
1 & 1 & \cdots & 1 & 1 \\ 
\alpha & \alpha & \cdots & \alpha & 1 \\
\alpha^2 & \alpha^2 & \cdots & \alpha & 1 \\ 
\vdots & \vdots & \vdots & \vdots & \vdots \\ 
\alpha^n & \alpha^{n-1} & \cdots & \alpha & 1 
\end{array} \right]$$ 
Then $Q_n(r,k) > \beta^k$ 
where $\beta= \beta_{r,n} = \sqrt{r/\lambda_{max}(A_{r,n})}$. Note that 
for each $r$, there are only finitely many values for which $\beta_{r,n} > 1$. 
The first few such values are shown in the following table.

\begin{center}

\begin{tabular}{|c|c|c|c|c|c|c|c|}
\hline
$n$ & $1$ & $2$ & $3$ & $4$ & $5$ & $6$ \\
\hline
$\beta_{2,n}$ & $1.08239$ & $< 1$ & $< 1$ & $< 1$ & $< 1$ & $< 1$ \\
\hline
$\beta_{3,n}$ & $1.28511$ & $1.11226$ & $1.02236$ & $< 1$ & $< 1$ & $< 1$ \\
\hline
$\beta_{4,n}$ & $1.46410$ & $1.24686$ & $1.12770$ & $1.05338$ & $1.00384$ 
& $< 1$ \\
\hline
\end{tabular}
\end{center}


\begin{thebibliography}{99}

\bibitem{Gow01} W. T. Gowers, {\em {A new proof of Szemer\'{e}di's
theorem}}. Geometric and Functional Analysis 11, 2001.

\bibitem{Job11} A. Jobson, A. Kezdy, H. Snevily and S. C. White, {\em 
{Ramsey functions for quasi-progressions with large diameter}}. Journal of 
Combinatorics 2 (2011), 557-573.

\bibitem{Lan98} B. M. Landman, {\em {Monochromatic sequences whose gaps 
belong to $\{d,2d,\ldots,md\}$}}. Bulletin of the Australian Mathematical 
Society 58 (1998), 93-101. 

\bibitem{Lan04} B. M. Landman and A. Robertson, {\em {Ramsey Theory on the 
Integers}}. American Mathematical Society, Providence, 2004.

\bibitem{Vij10} S. Vijay, {\em {On a variant of Van der Waerden's 
Theorem}}. Integers 10 (2010), A17, 5pp. (electronic).

\bibitem{Wae27} B. L. van der Waerden, {\em {Beweis einer Baudetschen
Vermutung}}. Niew Archief voor Wiskunde 15 (1927), 212-216.

\end{thebibliography}
\end{document}